\documentclass[12pt]{article}
\usepackage[centertags]{amsmath}
\usepackage{amsfonts}
\usepackage{amssymb}
\usepackage{amsthm}
\usepackage{amsmath,mathrsfs}
\usepackage{amsmath,amscd}
\usepackage[matrix,arrow,curve,cmtip]{xy}
\title{\textbf{Segal Topoi and Natural Phenomena: Universality of Physical Laws}}
\author{Renaud Gauthier \footnote{rg.mathematics@gmail.com} \\ \\}
\theoremstyle{definition}
\newtheorem{Thm2}{Theorem}[subsection]
\newtheorem*{acknowledgments}{Acknowledgments}

\DeclareMathOperator*{\colim}{\text{colim}}
\DeclareMathOperator*{\adj}{\rlarr}
\DeclareMathOperator*{\rel}{\sim}

\newcommand{\beq}{\begin{equation}}
\newcommand{\eeq}{\end{equation}}

\newcommand{\rarr}{\rightarrow}

\newcommand{\rlarr}{\rightleftarrows}

\newcommand{\xrarr}{\xrightarrow}

\newcommand{\cA}{\mathcal{A}}

\newcommand{\cC}{\mathcal{C}}
\newcommand{\cCop}{\cC^{\op}}

\newcommand{\cE}{\mathcal{E}}

\newcommand{\cG}{\mathcal{G}}
\newcommand{\cH}{\mathcal{H}}

\newcommand{\gC}{\mathfrak{C}}

\newcommand{\bL}{\mathbb{L}}

\newcommand{\bR}{\mathbb{R}}

\newcommand{\Cat}{\text{Cat}}

\newcommand{\diag}{\text{diag}}

\newcommand{\Hom}{\text{Hom}}
\newcommand{\Ho}{\text{Ho}\,}

\newcommand{\Map}{\text{Map}}

\newcommand{\op}{\text{op}}

\newcommand{\Set}{\text{Set}}

\newcommand{\Top}{\text{Top}}
\newcommand{\uHom}{\underline{\Hom}}

\newcommand{\AffC}{\cA \text{ff}_{\cC}}

\newcommand{\Catinf}{\Cat_{\infty}}

\newcommand{\CatD}{\Cat_{\Delta}}

\newcommand{\gCKop}{\gC[K]^{\op}}

\newcommand{\Kan}{\text{Kan}}

\newcommand{\Ktr}{K^{\triangleright}}

\newcommand{\LBous}{L_{\text{Bous}}}

\newcommand{\LSetDplusgCKop}{L(\SetD^+)^{\gCKop}}
\newcommand{\LgCKop}{L\gCKop}

\newcommand{\RHom}{\mathbb{R} \uHom}

\newcommand{\RuHom}{\bR \uHom}

\newcommand{\SetD}{\Set_{\Delta}}

\newcommand{\SePC}{\text{SePC}}

\newcommand{\SetDK}{(\SetD)_{/K}}
\newcommand{\SetDplusK}{(\SetD^+)_{/K}}
\newcommand{\SetDplusC}{(\SetD^+)^{\cC}}

\newcommand{\Stf}{\text{St}_f}
\newcommand{\Stplusf}{\Stf^+}

\newcommand{\SetDC}{(\SetD)^{\cC}}

\newcommand{\SePrStck}{\text{SePrStck}}

\newcommand{\tleft}{\triangleleft}

\newcommand{\Unf}{\text{Un}_f}
\newcommand{\Unplusf}{\Unf^+}

\newcommand{\Xtr}{X^{\triangleright}}

\begin{document}
\maketitle
\begin{abstract}
J. Lurie proved in \cite{Lu1} that for $K \in \SetD$, $\cC \in \CatD$, $f: \gC[K] \rarr \cC^{\op}$ an equivalence of simplicial categories, we have a Quillen equivalence $\Stplusf: \SetDplusK \rlarr \SetDplusC: \Unplusf$. We prove a partial converse to this theorem at the level of Segal categories, namely that if $L\SetDplusK$ is isomorphic to $L\SetDplusC$ in $\Ho(\SePC)$, then $\LgCKop$ and $L \cC$ are equivalent as Segal pre-categories relative to Segal categories of prestacks. We interpret this as indicating that the Segal category of pre-stacks $L \SetDplusC \cong \RuHom(LC, L \SetD^+)$ on $L\cC$ is equivalently given by a choice of simplicial set $K$, relative to which phenomena in $\Top^+ = L(\SetD^+)$ are considered, a sort of relativity principle. If we further take the Bousfield localizations of $L(\SetD^+)^{\gCKop} \cong L \SetDplusK$ and $L\SetDplusC$ with respect to hypercovers, then regarding $\LBous(L\SetDplusC)$ as the Segal topos of natural phenomena on $L\cC$, we also obtain an isomorphism $\LBous (L(\SetD^+)^{\gCKop}) \cong \LBous(L\SetDplusC)$ of Segal topoi of stacks. This provides two representations of the same natural phenomena, concurrently with the equivalence $\LgCKop \sim L\cC$ relative to prestacks, which we interpret as a weak universality of natural laws.
\end{abstract}

\newpage

\section{Introduction}
Lurie showed in \cite{Lu1} that starting from the usual Grothendieck construction whereby one obtains a category $\cC_f \xrightarrow{p} \cC$ cofibered in groupoids over $\cC$ from a functor $f: \cC \rarr \mathcal{G}rpd$ valued in groupoids, one gets, via the simplicial nerve functor, a left fibration in $\SetD$: $N(\cC_f) \xrightarrow{Np} N(\cC)$. Generalizing this to the $\infty$-categorical setting, Lurie was led to proving that one has a Quillen adjunction:
\beq
\SetDK \rlarr \SetDC \nonumber
\eeq
for $K \in \SetD$, $\cC \in \CatD$, $\gC[K] \rarr \cCop$ a morphism in $\CatD$. Here $\SetDK$ is endowed with the contravariant model structure (see Section \ref{Sec2}), and the projective model structure is put on $\SetDC$. Note that in the contravariant model structure it is right fibrations that correspond to fibrant objects. The choice of such a model structure is owing to the fact that right fibrations are technically easier to work with than are left fibrations. Further if $\gC[K] \rarr \cCop$ is an equivalence of simplicial categories, then one obtains a Quillen equivalence above.\\

Since our aim is to apply this formalism to fibered categories, it is preferable as pointed out in \cite{Lu1} to work with marked simplices instead. The bridge between both formalisms is warranted by virtue of the fact proved in \cite{Lu1} that we have a Quillen equivalence:
\beq
\SetDK \rlarr \SetDplusK \nonumber
\eeq
for all $K \in \SetD$, where we have put the contravariant model structure on $\SetDK$ and the Cartesian model structure on $\SetDplusK$ (see Section \ref{Sec4}). In the marked setting we have a Quillen adjunction:
\beq
\SetDplusK \rlarr \SetDplusC \nonumber
\eeq
for every morphism $f:\gC[K] \rarr \cCop$ in $\CatD$. Further if $f$ is an equivalence, this adjunction becomes a Quillen equivalence. 

\newpage

We prove a partial converse to this result, namely that after taking a simplicial localization of such simplicial categories (\cite{DK1}, \cite{TV1}), we have $L\SetDplusK$ is isomorphic to $L(\SetD^+)^{\gC[K]^{\op}}$, and that if the latter Segal category is isomorphic to $L\SetDplusC$ (both isomorphisms in $\Ho(\SePC)$, the homotopy category of Segal pre-categories, i.e. the category of Segal categories), then the morphism $\LgCKop \rarr L\cC$ induced by $\gC[K]^{\op} \rarr \cC$ is an equivalence in $\SePC$ relative to Segal categories of prestacks. \\

To go a little further, it is not a Segal topos of prestacks one should consider if one is interested in modeling natural phenomena. Rather, as argued in \cite{RG} and \cite{RG2}, it is really stacks we should work with. If we endow $L\gC[K]^{\op}$ with a Segal topology, and we consider morphisms of Segal sites $L\gC[K]^{\op} \rarr L \cC$, taking Bousfield localizations of $L\SetDplusC$ and $L(\SetD^+)^{\gC[K]^{\op}}$ with respect to their respective hypercovers, we get isomorphic Segal topoi of stacks. In other terms $\LSetDplusgCKop \cong L\SetDplusC$ implies getting isomorphic representations of natural phenomena after Bousfield localization, and the base Segal categories $L \gCKop$ and $L \cC$ giving rise to such phenomena are equivalent relative to prestacks, i.e. correspond to the same intrinsic natural laws. This would point to a universality of physical laws. Further this is for $K$ fixed, hence we consider phenomena in reference to $K$, via the isomorphism $L\SetDplusK \cong \LSetDplusgCKop$. Had we picked another simplicial set, we would have obtained another isomorphism $\LSetDplusgCKop \cong L\SetDplusC$, hence another homotopy class of $L\gC[K] \in \SePC$. In other terms we have a sort of relativity principle whereby the point of view (being $K$) dictates what are the physical laws that one should consider. This is partly the vision Grothendieck had as far back as \cite{G1}, namely that Physics should occur in higher categories.\\

The foundational results, and statements pertaining to the adjunctions whose converse we prove, are all covered in full in \cite{Lu1} and are just included here for ease of reading. No originality is claimed. We just reorganize the presentation slightly for our purposes, and notations may vary only slightly. Our references for model category theory are standard (\cite{Ho}, \cite{Hi}). For Segal categories we use \cite{TV1}, \cite{TV2}, \cite{T}, \cite{P} and \cite{HS}. For stacks we use \cite{TV} and \cite{TV4}.\\

\newpage

In Section \ref{Sec2} we discuss the contravariant model structure on $\SetDK$. In Section \ref{Sec3}, we give the unmarked Quillen equivalence $\SetDK \rlarr \SetDC$ of \cite{Lu1}. Before making the transition to the marked case, we introduce marked simplicial sets in Section \ref{Sec4} and present the cartesian model structure on $\SetDplusK$ for $K \in \SetD$ as done in \cite{Lu1}. In Section \ref{Sec5} we give the marked Quillen equivalence of \cite{Lu1}: $\SetDplusK \rlarr \SetDplusC$. In Section \ref{Sec6} we have to define $\Top^+ = L \SetD^+$ since we will be using the strictification theorem $L \SetDplusC \cong \RHom(L\cC, L \SetD^+)$. Finally in Section \ref{Sec7} we give the main results of the paper, everything before, with the exception perhaps of Section \ref{Sec6}, being initially covered elsewhere.\\

\begin{acknowledgments}
	The author would like to thank T. Pilling for ongoing discussions that are both stimulating and illuminating, as well as the referee of this paper who pointed out an inconsistent claim in a first version of this work.
\end{acknowledgments}

\section{Contravariant model structure on $\SetDK$} \label{Sec2}
For $K \in \SetD$ fixed, one puts a model structure on $\SetDK$ that Lurie refers to as the contravariant model structure. In order to do so, we need a few definitions. Following \cite{Lu1}, one defines the homotopy category of spaces by $\cH = W^{-1} \SetD$ where $W$ is the set of weak equivalences in $\SetD$. If $\cC \cG$ denotes the category of compactly generated, Hausdorff topological spaces, one has a functor $[\;]: \cC \cG \rarr \cH, X \mapsto [X]$, $[X]$ the homotopy class of $X$. One also has the geometric realization functor $|\;|: \SetD \rarr \cC \cG$. Denote the composition by $\gamma = [\;] \circ |\;|: \SetD \rarr \cH$. Applying $\gamma$ to each of the morphism spaces of a simplicial category $\cC$ gives a $\cH$-enriched category that one denotes by $h\cC$. For $K \in \SetD$, one defines the homotopy category $hK$ of $K$ by $hK = h\gC[K]$, where $\gC[K] \in \CatD$ is defined in \cite{Lu1}(we won't need its construction). Recall that if one endows $\SetD$ with the Joyal model structure, with monomorphisms as cofibrations, and weak equivalences those maps $K \rarr K'$ in $\SetD$ inducing equivalences $\mathfrak{C}[K] \rarr \mathfrak{C}[K']$ in $\CatD$, one has a Quillen equivalence $\gC: \SetD \rlarr \CatD: N$, where $N$ is the simplicial nerve functor (\cite{C}). One says $f: K \rarr L$ in $\SetD$ is a \textbf{categorical equivalence} if $hf: hK \rarr hL$ is an equivalence of $\cH$-enriched categories. Recall the join of two simplicial sets $K$ and $L$: if $I$ is a nonempty, finite, linearly ordered set, one has:
\beq
K \star L = \coprod_{I = J \cup J'} K(J) \times L(J') \nonumber
\eeq

In particular one defines the left cone $K^{\tleft}$ as being $\Delta^0 \star K$ (and the right cone by $\Ktr = K \star \Delta^0$). Now a map $X \rarr Y$ in $\SetDK$ is said to be a \textbf{covariant equivalence} if it induces a categorical equivalence:
\beq
X^{\tleft}\coprod_X K \rarr Y^{\tleft} \coprod_Y K \nonumber
\eeq
One says a map in $\SetDK$ is a \textbf{covariant cofibration} if it is a monomorphism of simplicial sets, and it is a \textbf{covariant fibration} if it has the right lifting property with respect to weak covariant cofibrations. Covariant fibrations, cofibrations and equivalences determine a left proper, combinatorial model structure on $\SetDK$ called the \textbf{covariant model structure} (\cite{Lu1}). The \textbf{contravariant model structure} on $\SetDK$ is obtained by considering the same cofibrations, contravariant equivalences $f$ are those maps such that $f^{\op}$ is a covariant equivalence in $(\SetD)_{/K^{op}}$, and contravariant fibrations are those maps $f$ such that $f^{\op}$ is a covariant fibration in $(\SetD)_{/K^{\op}}$.\\

\section{Unmarked Quillen equivalence} \label{Sec3}
The construction that follows is done in full in \cite{Lu1}, and is reproduced here for completeness' sake, albeit with different notations. Let $K \in \SetD$, $\cC \in \CatD$, $f: \gC[K] \rarr \cCop$. Let $*$ be the cone point of $\Xtr$, $X \in \SetDK$. By definition of $X$, one has an induced map $\gC[X] \rarr \gC[K]$. Consider the pushout diagram:
\beq
\xymatrix{
\gC[X] \ar[d] \ar[r] &\gC[K] \ar[r]^f &\cCop \ar[d] \\
\gC[\Xtr] \ar[rr] &&\Gamma
} \nonumber
\eeq
where:
\beq
\Gamma = \Gamma[K,\cC,f,X] = \gC[\Xtr] \coprod_{\gC[X]} \cCop \nonumber
\eeq
which can be seen as a correspondence:
\beq
\Gamma[K, \cC, f, X] = \text{Corr}[K,f,X](\cCop \rarr \{*\}) \nonumber
\eeq
which itself can be viewed as a functor:
\begin{align}
St_f X:\, & \cC \rarr \SetD \nonumber \\
&C \mapsto \Map_{\Gamma}(c,*) \nonumber
\end{align}
an element of $\SetDC$, for each $X \in \SetDK$, giving rise to a \textbf{straightening functor}:
\begin{align}
St_f: \SetDK & \rarr \SetDC \nonumber \\
X &\mapsto St_f X \nonumber
\end{align}
There exists an adjunction:
\beq
St_f: \SetDK \rlarr \SetDC: Un_f \nonumber
\eeq
which turns out to be a Quillen adjunction if one endows $\SetDK$ with the contravariant model structure and $\SetDC$ with the projective model structure. Further if $f$ is an equivalence of simplicial categories, then $St_f \dashv Un_f$ is a Quillen equivalence. All this is proved in \cite{Lu1}.\\

\section{Marked simplicial sets} \label{Sec4}
As observed in \cite{Lu1}, if one wants to regard objects $X \in \SetDK$ functorially, one is led to regarding $X \rarr K$ as a Cartesian fibration, fibrant object in $\SetDK$ if this category is endowed with a model structure other than the contravariant model structure, one which necessitates considering marked simplicial sets, which are essentially simplicial sets with a collection of marked edges. Given the interpretation we want to draw from $\SetDK$, we are led to considering this new model structure, called the \textbf{Cartesian model structure}, which builds up on the one we introduced prior.\\

\newpage

As proved in \cite{Lu1}, we have a Quillen equivalence:
\beq
\SetDK \rlarr \SetDplusK \nonumber
\eeq
for all $K \in \SetD$, as discussed in the introduction. We can see this result as a first step towards generalizing the work in the unmarked case to the marked one.\\

For $K, S \in \SetD$, $p:S \rarr K$ a map in $\SetD$, there exists some $K_{/p} \in \SetD$ defined by:
\beq
\Hom_{\SetD}(L,K_{/p}) = \Hom(L \star S, K) \nonumber
\eeq
for all $L \in \SetD$, where on the right hand side we consider only those maps $L \star S \rarr K$ in $\SetD$ whose restriction to $S$ is just $p$. For $\cC \in \Catinf$, if $p:\Delta^n \rarr \cC$ classifies an $n$-simplex $\sigma$ of $\cC$, then write $\cC_{/\sigma}$ for $\cC_{/p}$. In particular, if $p=[0] = \Delta^0 \rarr \cC$ is a point $X$ of $\cC$, then write $\cC_{/p} = \cC_{/X}$ and if $p:\Delta^1 \rarr \cC$ is a morphism $f: x \rarr y$ in $\cC$, then write $\cC_{/p} = \cC_{/f}$.\\

With those notations, if $p:X \rarr K$ is an inner fibration of simplicial sets, an edge $u:x \rarr y$ of $X$ is said to be \textbf{p-cartesian} if the induced map:
\beq
X_{/u} \rarr X_{/y} \times_{K_{/p(y)}} K_{/p(u)} \nonumber
\eeq
is a trivial Kan fibration. \\

Now one says a map $p:X \rarr K$ in $\SetD$ is a \textbf{Cartesian fibration} if it is an inner fibration, and if in addition for all $u: x \rarr y$ in $K$, for any vertex $y_0$ of $X$ such that $p(y_0) = y$, there is a $p$-cartesian edge $u_0: x_0 \rarr y_0$ of $X$ such that $p(u_0) = u$.\\

One defines a \textbf{marked simplicial set} to be a simplicial set $X$ along with a set $\cE$ of edges of $X$ (called marked edges) which contains all degenerate edges. A map $f:(X, \cE) \rarr (X', \cE')$ of marked edges is just a map $f: X \rarr X'$ of simplicial sets such that $f(\cE) \subseteq \cE'$. One denotes the category of marked simplicial sets by $\SetD^+$. If $K \in \SetD$, we denote by $K^{\forall}$ the marked simplicial set $K$ for which all its edges are marked. In particular for $K \in \SetD$, $(\SetD^+)_{/K^{\forall}}$ is denoted $\SetDplusK$. If $p:X \rarr K$ is a Cartesian fibration in $\SetD$, we denote by $X^{cart}$ the marked simplicial set $X$ where we have kept as marked edges its $p$-cartesian edges only.\\

Now $\SetD^+$ is cartesian closed, so for all $X,Y \in \SetD^+$ there is an internal object $Y^X$. We denote by $\underline{Y^X}$ the underlying simplicial set of $Y^X$, and by $\Map^{\forall}(X,Y)$ the marked simplicial set $(\underline{Y^X})^{\forall}$. As proved in \cite{Lu1}, we have that for $p:X \rarr Y \in \SetDplusK$, the following two conditions are equivalent:
for any cartesian fibration $L \rarr K$, the induced map:
\beq
\underline{(L^{cart})^Y} \rarr \underline{(L^{cart})^X} \nonumber
\eeq
is an equivalence in $\Catinf$, and the induced map:
\beq
\Map^{\forall}(Y,L^{cart}) \rarr \Map^{\forall}(X,L^{cart}) \nonumber
\eeq
is a homotopy equivalence in $\Kan$, the category of Kan complexes. Finally, a map $p:X \rarr Y$ in $\SetDplusK$ satisfying those equivalent conditions will be referred to as a \textbf{Cartesian equivalence}.\\

There is a left proper, combinatorial model structure on $\SetDplusK$ called the \textbf{Cartesian model structure} (\cite{Lu1}), whose cofibrations are those maps in $\SetDplusK$ whose underlying maps on the underlying simplicial sets are cofibrations, weak equivalences are Cartesian equivalences, and fibrations are as usual those maps that have the right lifting property with respect to weak cofibrations. Observe, as this was our aim, that $X \in \SetDplusK$ with this model structure is fibrant if and only if it is equivalent to some $L^{cart}$ for $L \rarr K$ a Cartesian fibration. Note that $\SetDplusK$ is a simplicial model category.\\

\section{Marked Quillen equivalence} \label{Sec5}
We have to introduce a bit of notation first. For $K \in \SetD$, $X \in \SetDK$, $f: \gC[K] \rarr \cCop$,  $\sigma: \Delta^n \rarr \Map_{\cCop}(C,D)$, denote by $\sigma^*$ the following induced map:
\beq
St_f X(D)_n = \Map_{\Gamma}(D,*)_n \rarr \Map_{\Gamma}(C,*)_n =  St_f X(C)_n \nonumber
\eeq
where:
\beq
\Gamma = \gC[\Xtr] \coprod_{\gC[X]} \cCop \nonumber
\eeq
as introduced in Section \ref{Sec3}. 

\newpage

Indeed, if we have $C \rarr D \in \cCop$, this induces $\Map_{\Gamma}(C,*) \leftarrow \Map_{\Gamma}(D,*)$. In particular for $\sigma \in \Map_{\cCop}(C,D)_n$, one has the induced map $\sigma^*$ as defined above.\\

Now let $c:\Delta^0 \rarr X$, for $X \in \SetDK$. Then consider:
\beq
\xymatrix{
\Delta^0 \ar@{.>}[drrr] \ar[r]^c &X \ar@{~>}[d] \ar[r] &K \ar@{~>}[d] \nonumber \\
&\gC[X] \ar[r] &\gC[K] \ar[r]_f &\cCop \nonumber
}
\eeq
Let $C$ be the object of $\cC$ thus obtained. One can extend this map as follows:
\beq
\xymatrix{
c\star id_{\Delta^0}: \Delta^1 \ar@{.>}[drr] \ar[r] & \Xtr \ar@{~>}[d] \nonumber \\
&\gC[\Xtr] \ar[r] &\Gamma = \{\cCop \rarr * \}
}
\eeq
where we have morally identified $\Gamma$ with the set of maps $\cCop \rarr *$ to emphasize that $c \star id$ corresponds to a map $C \rarr *$, element of $St_f X(C)$, which we denote by $\tilde{c}$.\\

Now consider an edge $u: c\rarr d$ of $X$ as in:
\beq
\xymatrix{
\Delta^1 \ar@{.>}[drrr]^U \ar[r]^u &X \ar@{~>}[d] \ar[r] &K \ar@{~>}[d] \nonumber \\
&\gC[X] \ar[r] &\gC[K] \ar[r]_f & \cCop
}
\eeq
giving rise to $U:C \rarr D$ in $\cCop$. One can extend this to:
\beq
\xymatrix{
u\star id_{\Delta^0}: \Delta^2 \ar@{.>}[drr] \ar[r] &\Xtr \ar@{~>}[d] \nonumber \\
 &\gC[\Xtr] \ar[r] &\Gamma \nonumber \\
}
\eeq
where the dotted map corresponds to:
\beq
\xymatrix{
C \ar[dr]_{\tilde{c}} \ar[rr]^U &&D \ar[dl]^{\tilde{d}} \nonumber \\
& \ast
}
\eeq
along with a map $\tilde{u}: \tilde{c} \rarr \tilde{d} \circ U = U^* \tilde{d}$.\\

Now for $u:d \rarr e$ an edge of $X$, one has a corresponding diagram as above:
\beq
\xymatrix{
D \ar[dr]_{\tilde{d}} \ar[rr]^U &&E \ar[dl]^{\tilde{e}} \nonumber \\
&\ast
}
\eeq
in $St_fX(D)$, along with a map $\tilde{u}: \tilde{d} \rarr \tilde{e} \circ U = U^* \tilde{e}$, i.e. $\tilde{u} \in St_f X(D)_1$. For $\sigma: \Delta^1 \rarr \Map_{\cCop}(C,D)$, one has an induced map $\sigma^*: St_f X(D)_1 \rarr St_f X(C)_1$. Define $\cE_f(C)$ to be the collection of  $\sigma^* St_fX(D)_1$, where it is understood here that we take the pullbacks of all those $\tilde{u}$ in $St_f X(D)_1$ as defined above, i.e. this is for all edges $u:d \rarr e$ in $X$, and this for all $\sigma$, for all $D's$, i.e. all edges $u:c \rarr d$ of $X$. If we limit ourselves to edges in $\cE$, this leads us to defining:
\begin{align}
St_f^+(X,\cE): \cC & \rarr \SetD^+ \nonumber \\
C & \mapsto (St_fX(C), \cE_f(C)) \nonumber
\end{align}
There exists a right adjoint $Un_f ^+$ to $St_f^+$:
\beq
St_f^+: \SetDplusK \rlarr \SetDplusC: Un_f^+ \label{mQe}
\eeq
a Quillen adjunction, where $\SetDplusK$ is endowed with the Cartesian model structure and $\SetDplusC$ with the projective model structure. This adjunction is furthermore a Quillen equivalence if $f: \gC[K] \rarr \cCop$ is an equivalence of simplicial categories (\cite{Lu1}).\\

\section{Simplicial localization of $\SetD^+$} \label{Sec6}
With a view towards taking a localization of the above Quillen equivalence \eqref{mQe}, we need a notion of localization of $\SetD^+$. In \cite{TV1}, \cite{T}, Toen and Vezzosi denote the simplicial localization (\cite{DK1}) $L\SetD$ of $\SetD$ by $\Top$. We define its marked counterpart. In order to do so we start from \cite{DK1}, and we adapt their construction to $\SetD^+$.\\

The free category on $\SetD^+$ is the category $F\SetD^+ \in \CatD$ with a generator $F \phi$ for every non-identity map $\phi \in \SetD^+$. This construction comes with two functors $D:F\SetD^+ \rarr \SetD^+$ and $U:F\SetD^+ \rarr F^2 \SetD^+$ satisfying the usual comonad conditions. With this in hand we can define the standard free resolution of $\SetD^+$ to be the simplicial object $F_* \SetD^+$ with $F_k \SetD^+ = F^{k+1} \SetD^+ \in \CatD$ with face and degeneracy operators given by:
\beq
d_i: F^{k+1} \SetD^+ \xrightarrow{F^iDF^{k-i}} F^k \SetD^+ \nonumber
\eeq
and:
\beq
s_i: F^{k+1} \SetD^+ \xrightarrow{F^i U F^{k-i}} F^{k+2} \SetD^+ \nonumber
\eeq
Observe that we have a weak equivalence of bi-simplicial categories:
\beq
F_* \SetD^+ \rarr \SetD^+ \nonumber
\eeq
given by:
\beq
F^{k+1} \SetD^+ \xrightarrow{D^{k+1}} \SetD^+ \nonumber
\eeq
To localize $ F_* \SetD^+$, we need a notion of weak equivalence on $\SetD^+$. Recall from \cite{GoJa} the definition of the simplicial category $\Delta \downarrow K$ of a simplicial set $K$: its objects are maps $\sigma: \Delta^n \rarr K$ with arrows being commutative diagrams:
\beq
\xymatrix{
\Delta^n \ar[dd]_{\Xi} \ar[dr]^{\sigma} \\
&K \\
\Delta^m \ar[ur]_{\tau}
} \label{SC}
\eeq
which gives an isomorphism $K \cong \colim \Delta^n$ where the colimit is taken over all morphisms $\Delta^n \rarr K$ in $\Delta \downarrow K$. This enables one to define the geometric realization of $K$ as $|K| = \colim |\Delta^n| \in \Top$ where the colimit is over the same maps $\Delta^n \rarr K$. Now for $K^+ = (K, \cE) \in \SetD^+$, using the same morphisms \eqref{SC}, one has $\sigma (\Delta^n) \cap \cE = \tau \circ \Xi (\Delta^n) \cap \cE$, i.e. $\cE|_{\sigma \Delta^n} = \cE|_{\tau \circ \Xi\Delta^n}$, hence the right vertical map in the diagram below is a morphism of marked simplicial sets:
\beq
\xymatrix{
\Delta^n \ar[d] \ar[r]^{\sigma} &(K, \cE) \ar[d]^{id_K} \nonumber\\
\Delta^m \ar[r]_{\tau} &(K,\cE)
}
\eeq
It seems one could define the simplicial category $\Delta \downarrow K^+$ as in the unmarked case. However we are really only keeping those simplices $\sigma: \Delta^n \rarr K$ who have an edge in $\cE$, i.e. those $n$-simplices in $K_n \times_{K_1} \cE$. 
This allows us to define the simplicial category $\Delta \downarrow K^+$ as the category whose objects are maps $\sigma: \Delta^n \times \Delta^0 \rarr (K,\cE)$, elements of $K_n \times_{K_1} \cE$, and morphisms are commutative diagrams:
\beq
\xymatrix{
\Delta^n \times \Delta^0 \ar[dd]_{\Xi} \ar[dr]^{\sigma} \nonumber \\
&(K,\cE) \nonumber \\
\Delta^m \times \Delta^0 \ar[ur]_{\tau} \nonumber
}
\eeq
which allows us to define the geometric realization of $K^+$ as:
\beq
|(K,\cE)| = \colim |\Delta^n \times \Delta^0| \nonumber
\eeq
where the colimit is taken in $\SetDplusK$, over those $\sigma : \Delta^n \times \Delta^0 \rarr K$ in $K_n \times_{K_1} \cE$. Weak equivalences in $\SetD^+$ between marked simplicial sets are then those maps that induce weak homotopy equivalences between their respective realizations. Let $W^+$ be the set of weak equivalences in $\SetD^+$. With this notion of equivalence we have a functor:
\beq
F_* \SetD^+ \rarr F_* \SetD^+[(F_*W^+)^{-1}] \nonumber
\eeq
defined levelwise. We then define:
\beq
\Top^+ = L(\SetD^+) = \diag F_* \SetD^+[(F_* W^+)^{-1}]  \in \CatD \nonumber
\eeq
following the localization of simplicial categories as done in \cite{DK1}.

\newpage

\section{Equivalence of models for natural phenomena} \label{Sec7}
Applying the Quillen equivalence \eqref{mQe} to the case where $\cC = \gC[K]^{\op}$ one obtains a Quillen equivalence:
\beq
St^+: \SetDplusK \rlarr (\SetD^+)^{\gC[K]^{\op}} :Un^+\nonumber
\eeq
In \cite{Lu1} it was proved that if $\gC[K] \rarr \cCop$ is an equivalence in $\CatD$, then one has a Quillen equivalence:
\beq
(\SetD^+)^{\gCKop} \rlarr \SetDplusC \nonumber
\eeq
We wish to prove the converse, namely that if we have such an equivalence, then the two simplicial categories we started from are equivalent in a sense to be precised. This we do after simplicial localization.

\subsection{Isomorphism of Segal categories of pre-stacks}
According to the strictification theorem (\cite{HS}, \cite{T}, \cite{TV}), we have an equivalence of Segal categories in $\SePC$:
\beq
L\SetDplusC \simeq \RHom(L\cC, \Top^+) \label{sim1}
\eeq
hence an isomorphism in $\Ho(\SePC)$. We use the isomorphism since Segal categories are objects of $\Ho(\SePC)$. Here $\RHom$ is the derived internal Hom in $\Ho(\SePC)$, with $\uHom$ the internal Hom in $\SePC$, the category of Segal pre-categories (see \cite{TV1} for notations). We will prove if $u: L \cC' \rarr L \cC$ induces an equivalence $L\SetDplusC \simeq L(\SetD^+)^{\cC'}$, then we have $L\cC \simeq L\cC'$ in $\SePC$ in a certain sense. For this we use the following result (\cite{Hi}): if $M \in \CatD$, $f: X \rarr Y$ a morphism in $M$ where both $X$ and $Y$ are cofibrant, then $f$ is a weak equivalence if and only if for any fibrant $Z \in M$ the induced map $f^*: \Map(Y,Z) \rarr \Map(X,Z)$ is a weak equivalence in $\SetD$. Here $\SePC$ is a simplicial model category to which we apply this result. Observe that morphisms of Segal categories are defined via the derived internal hom $\RHom $ in $\Ho(\SePC)$, while we denote by $\Map$ the mapping space object of $\SePC$. The equivalence $L\SetDplusC \simeq L(\SetD^+)^{\cC'}$ in $\SePC$, induces an isomorphism in $\Ho(\SePC)$. Further \eqref{sim1} above gives us another isomorphism in $\Ho(\SePC)$. Combining both:
\beq
\RuHom(L\cC, \Top^+) \cong \RuHom(L\cC', \Top^+) \nonumber
\eeq
which means that the Segal categories $L\cC$ and $L\cC'$ produce isomorphic Segal categories of prestacks in $\Ho(\SePC)$. We can say more. We have:
\beq
\xymatrix{
	\RHom(L\cC, \Top^+) \ar@{=}[d] \ar[r]^{\cong}_{\Ho(\SePC)} & \RHom(L\cC', \Top^+) \ar@{=}[d] \nonumber\\
\uHom(QL\cC, R\Top^+) \ar@{.>}[r] &\uHom(QL\cC', R\Top^+)
}
\eeq
where we have used the definition of the derived Hom $\RuHom$ and the fact that Segal categories are cofibrant in $\SePC$. Hence we have an equivalence in $\SePC$:
\beq
\uHom(QL\cC, R\Top^+) \simeq \uHom(QL\cC', R\Top^+) \nonumber
\eeq
which allows us to write, for $D \in \SePC$:
\beq
\xymatrix{
\Map(QD, \uHom(QL\cC, R\Top^+)) \ar[d]_{\simeq} \ar[r]^{\simeq} & \Map(QD, \uHom(QL\cC', R\Top^+)) \ar[d]^{\simeq} \nonumber \\
\Map(QD \times QL\cC, R\Top^+) \ar[d]_{\simeq} & \Map(QD \times QL\cC', R\Top^+) \ar[d]^{\simeq} \\
\Map(QL\cC, \uHom(QD, R\Top^+)) \ar@{.>}[r] &\Map(QL\cC', \uHom(QD, R\Top^+))
}
\eeq
Since all objects are cofibrant in $\SePC$, and $u: L\cC' \rarr L\cC$ induces the dotted map above, by the 2-3 property this reads:
\beq
\Map(L\cC, \hat{E}) \xrarr{\simeq} \Map(L\cC', \hat{E}) \nonumber
\eeq
for any prestack $\hat{E}$. This means $u: L\cC' \rarr L\cC$ is local with respect to Segal categories of prestacks. Note that this is different from the original definition of local object in \cite{Hi} and is just meant here to indicate that relative to Segal categories of prestacks, the Segal categories $L\cC$ and $L \cC'$ are considered to be equivalent. We denote this by 
\beq
L \cC' \rel_{\SePrStck} L \cC \nonumber
\eeq
Applying this to the case $\cC' = \gC[K]^{\op}$ in particular, we have $L(\SetD^+)^{\gC[K]^{\op}} \cong L\SetDplusC$ in $\Ho(\SePC)$ (induced by $\gC[K]^{\op} \rarr \cC$) implies an equivalence: 
\beq
\LgCKop \rel_{\SePrStck} L\cC \nonumber 
\eeq

\subsection{Statement of the theorem}
Since we have a Quillen equivalence $\SetDplusK \rlarr (\SetD^+)^{\gC[K]^{\op}}$, the hammock localizations of their subcategories of cofibrant objects are weakly equivalent (\cite{DK3}):
\beq
\diag L^H [(\SetD^+)_{/K}]^c \xrarr{\simeq} \diag L^H[(\SetD^+)^{\gCKop}]^c \nonumber
\eeq
Given the following equivalences:
\beq
\xymatrix{
	\diag L^H [(\SetD^+)_{/K}]^c \ar[d]_{\simeq} \ar[r]^{\simeq} & \diag L^H[(\SetD^+)^{\gCKop}]^c \ar[d]^{\simeq} \\
	\diag L^H (\SetD^+)_{/K} & \diag L^H(\SetD^+)^{\gCKop} 
} \nonumber
\eeq
it follows $\diag L^H (\SetD^+)_{/K} \xrarr{\cong} \diag L^H(\SetD^+)^{\gCKop}$. Further if $\cC \in \CatD$, then from \cite{DK1} we have a roof diagram of equivalences:
\beq
\xymatrix{
	& \diag L^H F_* \cC \ar[dl]_{\simeq} \ar[dr]^{\simeq} \\
	\diag L^H \cC && L\cC
} \nonumber 
\eeq
Collecting things, we have an isomorphism in $\Ho(\SePC)$:
\beq
L\SetDplusK \xrightarrow{\cong} L(\SetD^+)^{\gC[K]^{\op}}  \nonumber
\eeq
We have proved:
\begin{Thm2} \label{Thm2}
	If $L\SetDplusK \cong L\SetDplusC$ in $\Ho(\SePC)$, $u: \gC[K]^{\op} \rarr \cC$, then $\LgCKop \rel_{\SePrStck} L\cC$ in $\SePC$.
\end{Thm2}

\subsection{Determining $K$ in $L\SetDplusK$}
In this section, we are interested in the peripheral question of finding what reference simplicial set $K$ would $\cC \in \CatD$ correspond to when we consider the isomorphism $L(\SetD^+)_{/K} \cong L(\SetD^+)^{\cC}$. Recall from \cite{Lu1} that one has a Quillen adjunction $\gC: \SetD \rlarr \CatD: N$, hence for any $\cC \in \CatD$, a weak equivalence: $\bL \gC(N(R\cCop)) \rarr R\cCop$. Since every object of $\SetD$ is cofibrant, we have an equivalence $\gC[K]^{\op} \rarr R\cC$ in $\CatD$ with $K = N(R \cCop)$ (\cite{Ho}). For simplicity, we will consider $\cC$ fibrant in $\CatD$. Because we have such an equivalence, by \cite{Lu1} we have a Quillen equivalence $(\SetD^+)_{/K} \adj (\SetD^+)^{\cC}$. Following the exact same reasoning as in the previous section, it follows that we have an isomorphism of Segal categories $ L (\SetD^+)_{/K} \xrarr{\cong} L(\SetD^+)^{\cC}$ showing that $K = N(\cCop)$ is the desired simplicial set.\\

\subsection{Physical phenomena}
To push the result of Theorem~\ref{Thm2} further, physical phenomena occur in Segal topoi of stacks, not prestacks, so one has to localize those equivalences. One needs to first put a topology on the Segal categories $\LgCKop$ and $L\cC$. Recall from \cite{TV1} that a topology on $L\cC$ is a Grothendieck topology on its homotopy category $\Ho(L\cC)$. Let $\tau$ be one such topology. Denote $(L\cC)^{\op}$ by $\AffC$, and $\AffC^{\sim, \tau} = \LBous \RHom(L\cC, \Top^+)$, the Segal category of stacks on $\AffC$, where the Bousfield localization is with respect to hypercovers for the topology $\tau$. Now the local equivalence $\LgCKop \rel_{\SePrStck} L\cC$ in $\SePC$ induces a morphism of homotopy categories $\rho: \Ho(\LgCKop) \rarr \Ho(L\cC)$. We are interested in having such a map preserve the Segal categories of stacks on those sites. This happens if $\rho$ is a morphism of sites, or a continuous morphism as it is sometimes called. If that's the case, one can take simultaneous Bousfield localizations as in:
\beq
\xymatrix{
\RHom(\LgCKop, \Top^+) \ar[d]_{\LBous} & \RHom(L\cC, \Top^+) \ar[l]_-{\cong} \ar[d]^{\LBous} \nonumber \\
L\gC[K]^{\sim, \tau'} & \AffC^{\sim, \tau} \ar@{.>}[l]_{\cong} 
}
\eeq
From the top isomorphism then, we have isomorphic phenomena, modeled by Segal topoi of stacks on the bottom row, which correspond to equivalent Segal categories $\LgCKop$ and $L\cC$ relative to Segal categories of prestacks according to Theorem~\ref{Thm2}.\\

Note as an aside that this implies we do not have a unique set of natural laws, but a homotopical class of such laws.\\

\newpage

To summarize, we have an isomorphism of Segal topoi of stacks, concurrently with an equivalence of simplicial sites, both of which follow from an isomorphism in $\Ho(\SePC)$ of Segal categories of pre-stacks. We regard this formalism as a roof diagram:
\beq
\xymatrix{
&\LSetDplusgCKop \ar[ld]_{\LBous} \ar@{.>}[dr] &L \SetDplusC \ar[ld]^{\LBous} \ar[l]_-{\cong} \ar@{.>}[dr] \nonumber \\
(L \mathfrak{C}[K])^{\sim, \tau'} & \mathcal{A}ff_{\cC}^{\sim, \tau} \ar@{.>}[l]_-{\cong} &\LgCKop \ar[r]^{\sim}_{\SePrStck} & L\cC
}
\eeq
which we interpret as a weak universality of natural laws, the choice of word universality being meant to indicate that Segal sites in a same class yield equivalent physical phenomena.

\end{document}